\documentclass{article}
\usepackage[T1]{fontenc}
\usepackage{amsmath, amssymb}
\usepackage[hypertex]{hyperref}
\usepackage[authoryear,  comma, sort&compress]{natbib}
\RequirePackage{hypernat}
\usepackage{graphicx,   epic, eepicemu, epsfig, footmisc, float, enumerate}
\newtheorem{theorem}{Theorem}[section]

\newtheorem{lemma}{Lemma}[section]
\newtheorem{example}[theorem]{Example}

\newtheorem{corollary}{Corollary}[section]
\def\X{\mathcal{X}}

\newenvironment{proof}[1][Proof]{\textbf{#1.} }{\ \rule{0.5em}{0.5em}}

\def\title#1{{\Large\bf  \begin{center} #1 \vspace{3pt} \end{center}  } }
\def\authors#1{{\large\bf \begin{center} #1 \vspace{1pt} \end{center} } }
\def\university#1{{\sl \begin{center} #1 \vspace{5pt} \end{center} } }
\def\inst#1{\vspace{1pt} \unskip$^{#1}$}

\begin{document}\allowdisplaybreaks

%
%
%
\bibliographystyle{elsart-harv}

\title{Infinite Viterbi alignments in the two state hidden
Markov models}
\authors{J\"uri Lember\inst{1}\footnote{The author has been supported by
     the Estonian Science
 Foundation Grant 7553} and Alexey Koloydenko\inst{2} }
\university{\inst{1} University of Tartu, Estonia, email:jyril@ut.ee \\
 \inst{2} Royal Holloway, University of London, UK,  email:alexey.koloydenko@rhul.ac.uk}

\bigskip
\noindent {\large\bf Abstract}

\smallskip

\noindent Since the early days of digital communication, Hidden Markov Models
(HMMs) have now been routinely used in speech recognition, processing of
natural languages, images, and in bioinformatics.  An HMM
$(X_i,Y_i)_{i\ge 1}$ assumes observations $X_1,X_2,\ldots$ to be
conditionally independent given an ``explanotary'' Markov process
$Y_1,Y_2,\ldots$, which itself is not observed; moreover, the
conditional distribution of $X_i$ depends solely on $Y_i$. Central
to the theory and applications of HMM is the Viterbi algorithm to find
{\em a maximum a posteriori} estimate $q_{1:n}=(q_1,q_2,\ldots,q_n)$ of
$Y_{1:n}$ given the observed data $x_{1:n}$.
Maximum {\em a posteriori} paths are also called 
Viterbi paths or alignments.  Recently, attempts have been made to study the behavior of
Viterbi alignments of HMMs with two hidden states when $n$ tends to infinity. It has indeed been
shown that in some special cases a well-defined limiting Viterbi alignment exists. While
innovative, these attempts have relied on rather strong assumptions.  This work proves the
existence of infinite Viterbi alignments for virtually any HMM with
two hidden states.
\bigskip

\noindent {\large\bf Keywords}

\smallskip

\noindent Hidden Markov models, maximum a posterior path, Viterbi alignment, Viterbi extraction, Viterbi training

\bigskip

\section{Introduction}
We consider  hidden Markov models (HMM) $(Y,X)$ with two hidden states.
Namely, $Y$ represents the hidden process $Y_1, Y_2, \ldots,$ which is an irreducible aperiodic Markov chain  with
state space $S=\{a,b\}$. In particular, the transition probabilities $\mathbb{P}=(p_{lm})$, $l,m\in S$,
are positive and the stationary distribution $\pi=\pi\mathbb{P}$ is unique. 
For technical convenience, $Y_1$ is assumed to follow $\pi$,
however, the results of the paper hold for arbitrary initial distributions.
To every state $l\in S$ there corresponds an {\it emission} distribution
$P_l$ on $\mathcal{X}=\mathbb{R}^d$.  Given a realization $y_{1:\infty}\in S^{\infty}$ of $Y$, the observations
$X_{1:\infty}:=X_1,X_2,\ldots$ are generated as follows. If $Y_i=a$ (resp.
$Y_i=b$), then $X_i$ is distributed according to $P_a$ (resp. $P_b$) and
independently of everything else. We refer to this model
as the (general) 2-state HMM.

In \citep{HMM}, HMMs are called `one of the most successful
statistical modelling ideas that have [emerged] in the last forty
years'. Since their classical application to digital communication
in 1960s (see further references in \citep{HMM}), HMMs have had a
defining impact on the mainstream research in speech recognition
\citep{vanaraamat, jelinek0, jelinek, mitEM2,  philips, IBM2002,
raamat, Rabiner86, mitEM1, philips2,  strom},  natural language
models  \citep{bilmes2, ochney}, and  more recently computational
biology  \citep{BioHMM2, natHMM, findgene, BioHMM1}. Thus, for example,
DNA regions can be labeled as $a$, `coding', or $b$, `non-coding',
with $P_a$ and $P_b$ representing the respective distributions on the
$\{A,C,G,T\}$ alphabet.

Given observations $x_{1:n}:=x_1,\ldots,x_n$,
and treating the hidden states $y_{1:n}:=y_1,\ldots,y_n$ as
parameters,  inference in HMMs typically involves $v(x_{1:n})$, a {\it maximum a posteriori} (MAP)
estimate of $Y_{1:n}$. It has now been recognized that  `[in] spite of the theoretical and
practical importance of the MAP path estimator, very little is known about its
properties' \citep{caliebe2}.  The same estimates are also known as {\it Viterbi}, or {\em forced
alignments} and can be efficiently computed by a dynamic programming algorithm
also bearing the name of {\em Viterbi}.  When substituted  for true $y_{1:n}$ in the
likelihood function $\Lambda(y_{1:n}; x_{1:n},\psi)$,
Viterbi alignments can also be used to estimate  $\psi$,
any unknown free parameters of the model.
Starting with an initial guess $\psi^{(0)}$ and alternating between
maximization of the likelihood $\Lambda(y_{1:n}; x_{1:n},\psi)$ in
$y_{1:n}$ and $\psi$ is at the core of {\em Viterbi training} (VT),
or {\em extraction} \citep{jelinek0}, also known as {\em segmental
K-means} \citep{Rabiner90, HMP}. Resulting estimates $\hat
\psi_{\mathrm{VT}}(x_{1:n},\psi^{(0)})$  are known to be different
from the {\em maximum likelihood} (ML) estimates $\hat
\psi_{\mathrm{ML}}(x_{1:n},\psi^{(0)})$  which in this case are most
commonly delivered by the EM procedure \citep{baumhmm, HMP,
emtutorial}.  Even if $\psi$ were  known, Viterbi alignments
$v(x_{1:n};\psi)$ would  typically differ from true paths $y_{1:n}$,
and the long-run properties of $v(x_{1:n};\psi)$ are not obvious
\citep{caliebe1, caliebe2, AVT1, AVT4, AVT3}.  Furthermore,
\citep{AVT1, AVT4, AVT3} propose a hybrid of VT and EM which  takes
into account the asymptotic discrepancy between $\hat
\psi_{\mathrm{ML}}(x_{1:n},\psi^{(0)})$ and $\hat
\psi_{\mathrm{VT}}(x_{1:n},\psi^{(0)})$ in order to increase
computational and statistical efficiencies of estimation of $\psi$
for $n$ large.  Thus or otherwise, an important question is {\em how
to find the asymptotic properties of Viterbi alignments, given that
$(n+1)^{\mathrm{th}}$ observation can in principle change the 
previous alignment entirely, i.e. $v(x_{1:n+1})_i\neq v(x_{1:n})_i$,
$1\le i\le n$?} Do the Viterbi alignments then admit well-defined
extensions? We answer this question positively in \citep{AVT4} for
general HMMs (in particular, allowing more than two hidden states)
by constructing {\em proper infinite Viterbi alignments}.
Generalizing and clarifying related results of \citep{caliebe1,
caliebe2}, the  approach in \citep{AVT4} is to extend alignments
{\em piecewise}, separating individual pieces by {\em nodes} (see
\S\ref{sec:prelim} below). Although the construction is natural, a
detailed formal proof of its correctness for general HMMs is rather
long  and requires certain mild technical assumptions. This  paper, on the other hand,
shows that {\em in the special case of two state HMMs, the existence
of infinite Viterbi alignments needs no special  assumptions and can
be proven considerably more easily.} The results of this paper
essentially complete and generalize those of \citep{caliebe1,
caliebe2}.
\section{Preliminaries}\label{sec:prelim}
Let $\lambda$ be a suitable $\sigma$-finite reference measure on $\mathbb{R}^d$
so that $P_a$ and $P_b$ have densities with respect to $\lambda$.
For example, $\lambda$ can be a Lebesgue measure, or, as in the case
of discrete observations, a counting measure. Thus, let $f_a$ and $f_b$ be the
densities of $P_a$ and $P_b$, respectively.
Throughout the rest of the paper, we assume that
$P_a\ne P_b$ or, equivalently,
\begin{equation}\label{A3}
\lambda\{x\in \mathcal{X}: f_a(x)\ne f_b(x)\}>0.
\end{equation}
Assumption \eqref{A3} is natural since there would be no need to
model the observation process by an HMM should the emission
distributions coincide. Note also that unlike in the general case,
the positivity of the transition probabilities is also a natural
assumption  for the two state HMMs. {\em No more assumption on the
HMM is made in this paper.} In particular, unlike \citep{caliebe2,
caliebe1}, we do not assume the square integrability of
$\log(f_a/f_b)$, or equality of the  supports of $P_a$ and $P_b$.
However, the latter  condition is not very restrictive, since for
the two state HMMs with unequal supports the existence of infinite
Viterbi alignments follows rather trivially (Corollary \ref{strong}).

Thus, for any $n\ge 1$ and any $x_{1:n}\in \mathcal{X}^n$ and $y_{1:n}\in S^n$,  the likelihood
$\Lambda_\pi(y_{1:n}; x_{1:n})$ is given by
$$\mathbf{P}(Y_{1:n}=y_{1:n})\prod_{i=1}^nf_{y_i}(x_i),~\text{where}~
\mathbf{P}(Y_{1:n}=y_{1:n})=\pi_{y_1}\prod_{i=2}^np_{y_{i-1}y_i}.$$
 Since estimation of $\psi$ is not a goal of this 
paper, the dependence on $\psi$ is suppressed. 
Decomposition \eqref{eqn:decompose} and 
recursion \eqref{eqn:recurse} below provide a basis for the Viterbi algorithm to compute
alignments. Namely, for all $u\in\{1,2,\ldots,n-1\}$, 
\begin{eqnarray}\label{eqn:decompose}
\max_{y_{1:n}\in S^{n}}\Lambda_\pi(y_{1:n}; x_{1:n})=
\max_{l\in S}\left[\delta_{u}(l)\times \hspace*{-1mm}\max_{y_{u+1:n}\in S^{n-u}}\hspace*{-1mm}
\Lambda_{(p_{l\cdot})}(y_{u+1:n}; x_{u+1:n})\right],&&
\end{eqnarray}
where $(p_{l\cdot})$ is the transition distribution given state
$l\in S$, and the {\it scores}
$$\delta_{u}(l):=\max_{y_{1:u-1}\in S^{u-1}}\Lambda((y_{1:u-1},l); x_{1:u}),\quad l=a,b,$$
are defined  for all $u\ge 1$, and $x_{1:u}\in \mathcal{X}^u$.
Thus, $\delta_u(l)$ is the maximum of the likelihood
of the paths terminating at $u$ in state $l$. Note that $\delta_1(l)=\pi_lf_l(x_1)$ and
$\delta_u(l)$ depends on $x_{1:u}$.
\begin{eqnarray}\label{eqn:recurse}
\delta_{u+1}(a)&=&\max\{\delta_u(a)p_{aa},\delta_u(b)p_{ba}\}f_a(x_{u+1}), \\
\delta_{u+1}(b)&=&\max\{\delta_u(a)p_{ab},\delta_u(b)p_{bb}\}f_b(x_{u+1}),~u\ge 1, \nonumber
\end{eqnarray}
\begin{example}\label{example:mixtures} 
Let $X_1,X_2,\ldots$ be i.i.d. following a mixture distribution $\pi_aP_a+\pi_bP_b$ with  
density $\pi_a f_a(x;\theta_a)+\pi_b f_b(x;\theta_b)$ and mixing
weights $\pi_a, \pi_b>0$. Such a sequence is  an HMM with
the transition probabilities   $\pi_a=p_{aa}=p_{ba}$, $\pi_b=p_{bb}=p_{ab}$.
In this special case the alignment is easy
to exhibit. 
Indeed, in this case  recursion \eqref{eqn:recurse} writes for any $u\ge 1$  as 
\begin{eqnarray}\label{eqn:recursem}
\delta_{u+1}(a)=c\pi_{a}f_a(x_{u+1}), & & 
\delta_{u+1}(b)=c\pi_{b}f_b(x_{u+1}),
\end{eqnarray}
where $c=max\{\delta_{u}(a),\delta_{u}(b)\}$.
Hence, the alignment  
$v(x_{1:n})$ can be obtained pointwise as follows:
$$v(x_{1:n})=(v(x_1), \ldots, v(x_n)),
\text{where}~v(x)=\arg\max\{\pi_a f_a(x), \pi_b f_b(x)\}.$$ 
Equivalently (ignoring possible ties), using a generalized
{\it Voronoi partition} $\X=X_a\cup \X_b$ with 
\begin{eqnarray*}
\X_a=\{x\in\mathcal{X}: \pi_a f_a(x)\geq \pi_bf_b(x)\},& &
\X_b=\{x\in\mathcal{X}: \pi_b f_b(x)> \pi_af_a(x)\},
\end{eqnarray*}
$v(x)=a$ if and only if $x\in \X_a$, and otherwise (i.e. $x\in \X_b$) $v(x)=b$. 
\end{example}
Generally, it follows from \eqref{eqn:recurse} that, if
\begin{equation}\label{strongnode-a}
\delta_u(a)p_{aa}>\delta_u(b)p_{ba},\quad
\delta_u(a)p_{ab}>\delta_u(b)p_{bb},
\end{equation}
for some $u$, $1\le u$, and some $x_{1:u}\in \mathcal{X}^{u}$, 
then for any $n>u$ and for any extension $x_{u+1:n}\in \mathcal{X}^{n-u}$, the Viterbi
alignment goes through state $a$ at time $u$.  Namely, truncation $v(x_{1:n})_{1:u}$ 
coincides with the Viterbi alignment  $v(x_{1:u})$ (indeed, \eqref{strongnode-a}
implies $\delta_u(a)>\delta_u(b)$). 
Thus, under condition \eqref{strongnode-a}, maximization of $\Lambda_\pi((y_{1:n},l); x_{1:n})$ can 
be reset at time $u$ by clearing $x_{1:u}$ from the memory,
retaining $v_{1:u}$, and replacing the initial distribution $\pi$ by $(p_{a\cdot})$ for
further maximization of  $\Lambda_{(p_{a\cdot})}(y_{u+1:n}; x_{u+1:n})$.
Following \citep{AVT4}, if condition \eqref{strongnode-a} holds,
then $x_u$ is called a {\it strong $a$-node} (of realization $x_{1:n}$, $n>u$), 
where `strong' refers to the inequalities in \eqref{strongnode-a} being strict. 

Suppose  $x_{1:\infty}$ contains infinitely  many strong $a$-nodes at times
$u_1<u_2<\ldots$. 
Let $v^1=v(x_{1:u_1})$, and let $v^k$ maximize
$\Lambda_{(p_{a\cdot})}(y_{u_{k-1}+1:u_k}; x_{u_{k-1}+1:u_k})$, for
all $k\ge 2$.  Then, concatenation $(v^1, v^2, v^3, \ldots)$ is
naturally called the {\em infinite piecewise Viterbi alignment}
\citep{AVT4}. Apparently, the {\em almost sure} existence of our
infinite alignments directly dependends on the  existence of
infinitely many (strong) nodes. At the same time,  whether or not
$x_u$ is a node depends on $x_{1:u}$ and hence is difficult to verify
directly. Fortunately, in many cases $x_u$ is guaranteed to be a
node based on several preceding observations $x_{u-m:u}$, $1\le
m<u$, ignoring the rest. Specifically, suppose for example that
$x\in \mathcal{X}$ is such that
\begin{equation}\label{verystrong}
p_{ia}f_a(x)p_{aj}>p_{ib}f_b(x)p_{bj}, \quad \forall
i,j\in S.\end{equation} It is easy to check that for any $u\ge 2$,  $x_u=x$ is a
strong $a$-node for any $x_{1:u-1}$. Hence, if $x_{1:\infty}$ contains infinitely many observations
satisfying \eqref{verystrong}, then $x_{1:\infty}$ also contains infinitely many strong nodes. This previous
condition  in its turn is met provided
\begin{equation}\label{verystronglambda}
\lambda\left(\left\{x\in \mathcal{X}: p_{ia}f_a(x)p_{aj}>p_{ib}f_b(x)p_{bj},\quad \forall
i,j\in S\right\}\right)>0.
\end{equation}
Indeed, since our underlying Markov chain $Y$ is ergodic, it is
rather easy to see that $X$ is ergodic as well \citep{HMP, volatility, leroux}.
Also, \eqref{verystronglambda} implies that $$P_a\left(\left\{x\in \mathcal{X}:
p_{ia}f_a(x)p_{aj}>p_{ib}f_b(x)p_{bj},\quad \forall i,j\in S\right\}\right)>0.$$ Thus,
it follows  from ergodicity of $X$ that almost every realization of $X$ has infinitely
many elements satisfying (\ref{verystrong}) and, hence infinitely many
strong nodes. We have thus proved the following Lemma.
\begin{lemma}\label{lemma0}
Assume that \eqref{verystronglambda} holds. Then almost every sequence
of observations $x_{1:\infty}$ has infinitely many strong $a$-nodes.\end{lemma}
(Clearly, interchanging $a$ and $b$ gives the same results in terms
of $b$-nodes.) Lemma \ref{strong} is essentially Theorem 1 in
\citep{caliebe1} (disregarding a misprint in the statement).
Condition \eqref{verystronglambda} holds for many two-state HMMs
including the so-called additive Gaussian noise model
\citep{caliebe2}, where the emission distributions are Gaussian.
Another trivial example is the model with  unequal supports of $P_a$
and $P_b$. Indeed, in that case \eqref{verystronglambda} holds (at least 
up to swapping $a$ and $b$). Hence, the following Corollary.
\begin{corollary}\label{strong}
If the supports of $P_a$ and $P_b$ are not equal, then  almost every
sequence of observations has infinitely many strong
nodes.\end{corollary}
{\em The goal of this work is essentially to remove condition
\eqref{verystronglambda} from Lemma \ref{lemma0}.}

To this end, following \citep{AVT4}, we call an observation satisfying
(\ref{verystrong}) an $a$-{\it barrier} of length 1. More generally,
a block of observations $z_{1:k}\in \mathcal{X}^k$ is called a (strong) barrier of length
$k\ge 1$ if for every $m\ge 0$ and $x_{1:m}\in \mathcal{X}^m$, $z_{1:k}$ contains a
(strong) node of realization $(x_{1:m},z_{1:k})$. In \citep{AVT4}, we
prove  the existence of infinitely many
barriers for a very general class of HMMs. For the two-state HMMs, the
conditions of our result in \citep{AVT4} are given by \eqref{ab} and \eqref{ba} below.
\begin{align}\label{ab} & P_a\left(\left\{x\in \mathcal{X}:f_a(x)\max\{p_{aa},p_{ba}\}>
f_b(x)\max\{p_{bb},p_{ab}\}\right\}\right)>0\quad{\rm and}\\
\label{ba} & P_b\left(\left\{x \in \mathcal{X}: f_b(x)\max\{p_{bb},p_{ab}\}>
f_a(x)\max\{p_{aa},p_{ba}\}\right\}\right)>0.
\end{align}
To achieve our goal, we will first prove the same result  for the
two-state HMM under the relaxed assumption that  (\ref{ab}) {\em or}
(\ref{ba}) holds. As we shall see below (Lemma~\ref{basic}), in our two-state HMM one of
these conditions is automatically satisfied and, moreover, all
barriers are strong. Hence, occurrence of infinitely many strong
barriers in this case will be shown (Theorem~\ref{maintheorem}) to require no additional assumptions.

Finally, if a node is not strong and $v(x_{1:n})$ is not unique,
an alignment might exist that does not go through this node. 
Such type of pathologies cause technical inconveniences in defining an infinite Viterbi alignment
and are treated in \citep{AVT4}. Fortunately, unlike in the general case, in the case  of two-state 
HMMs almost every realization has infinitely many strong nodes (Theorem~\ref{maintheorem}).  This
allows for a simple resolution of the non-uniqueness  in the case  of two-state 
HMMs. 
\section{Main results}
\subsection{Three types of the two-state HMM}
The following three cases exhaust all  the possibilities:
\begin{enumerate}
  \item $p_{aa}>p_{ba}\quad (\Leftrightarrow\quad p_{bb}>p_{ab})$;
  \item $p_{aa}<p_{ba}\quad (\Leftrightarrow\quad p_{bb}<p_{ab})$;
  \item $p_{aa}=p_{ba}\quad (\Leftrightarrow\quad p_{bb}=p_{ab})$.
\end{enumerate}
From the definition of nodes, it follows that  $x_u$ is {\em not} a
node only in one of the following two cases:
\begin{align*}
(A)\left\{%
\begin{array}{ll}
    \delta_u(a)p_{aa}> \delta_u(b)p_{ba} \\
    \delta_u(b)p_{bb}>\delta_u(a)p_{ab}
\end{array}%
\right.\quad {\rm or}\quad
(B)\left\{%
\begin{array}{ll}
    \delta_u(b)p_{ba}>\delta_u(a)p_{aa} \\
    \delta_u(a)p_{ab}>\delta_u(b)p_{bb}
\end{array}%
\right.
\end{align*}
Case (A) is equivalent to
\begin{equation}\label{yks}{p_{bb}\over p_{ab}}>{\delta_u(a)\over \delta_u(b)}>{p_{ba}\over p_{aa}}\end{equation}
and case (B) is equivalent to
\begin{equation}\label{kaks}
{p_{bb}\over p_{ab}}<{\delta_u(a)\over \delta_u(b)}<{p_{ba}\over
p_{aa}}.\end{equation} Thus, in case (A), we have
$\delta_{u+1}(a)=\delta_u(a)p_{aa}f_a(x_{u+1})$ and
$\delta_{u+1}(b)=\delta_u(b)p_{bb}f_b(x_{u+1})$, so that for any
$n>u$, the Viterbi alignment $v(x_{1:n})$ must satisfy
$v(x_{1:n})_{u}=v(x_{1:n})_{u+1}$. Similarly, in case (B)
$\delta_{u+1}(a)=\delta_u(b)p_{ba}f_a(x_{u+1})$ and
$\delta_{u+1}(b)=\delta_u(a)p_{ab}f_b(x_{u+1}),$ i.e.
$v(x_{1:n})_{u}\neq v(x_{1:n})_{u+1}$. Evidently, case 1 and case
(B) are mutually exclusive, and so are case 2 and case (A).
Therefore, if the transition matrix satisfies the conditions of case
1, then $x_u$ is not a node if and only if conditions (A) are
fulfilled. This implies that  in case 1, nodes are the only
possibility for $v(x_{1:n})$ to change state. On the other hand, if
the transition matrix satisfies the conditions of case 2, then $x_u$
is not a node if and only if (B) holds. Hence, in case 2 nodes are
the only possibility for $v(x_{1:n})$ to remain in one state. 
Case 3 corresponds to the mixture model (see Example \ref{example:mixtures} above).
Apparently \eqref{eqn:recursem}, every observation is a node in this case (see also 
Figure \ref{fig:2state}  below).

Let us now examine conditions \eqref{ab} and \eqref{ba}. From
equation (\ref{A3}), it follows that
\begin{equation}\label{A3equal} \lambda\left(\left \{x\in \mathcal{X}: f_a(x)>f_b(x)\right\}\right)>0,
\quad \lambda\left(\left \{x\in \mathcal{X}: f_a(x)<f_b(x)\right\}\right)>0
\end{equation}
and, for  any $\alpha>\beta>0$,
\begin{eqnarray}\label{A3unequal-a}
\lambda\left(\left\{x\in \mathcal{X}: \alpha f_a(x) > \beta f_b(x)\right\}\right)>0
\Leftrightarrow  P_a\left(\left\{x\in \mathcal{X}: \alpha f_a(x) > \beta
f_b(x)\right\}\right)>0& & \\\label{A3unequal-b} \lambda\left(\left\{x\in \mathcal{X}: \alpha f_b(x)>
\beta f_b(x)\right\}\right)>0 \Leftrightarrow  P_b\left(\left\{x\in \mathcal{X}: \alpha
f_b(y) > \beta f_b(y)\right\}\right)>0.&&
\end{eqnarray}
Therefore, we have the following Lemma.
\begin{lemma}\label{basic}
Any two state HMM satisfies at least one of the condtions \eqref{ab} and \eqref{ba}.
\end{lemma}
\begin{proof}
\underline{In case 1}, (\ref{ab}) and (\ref{ba}) are equivalent to
\begin{eqnarray}\label{eng}
P_a\left(\left\{x\in \mathcal{X}:f_a(x)p_{aa}>
f_b(x)p_{bb}\right\}\right)=P_a\left(\left\{x\in \mathcal{X}:{f_b(x)p_{bb}\over f_a(x)p_{aa}}< 1 \right\}\right)>0&&\\
\label{ing} P_b\left(\left\{x\in \mathcal{X}:f_b(x)p_{bb}>
f_a(x)p_{aa}\right\}\right)=P_b\left(\left\{x\in \mathcal{X}:{f_a(x)p_{aa}\over f_b(x)p_{bb}}< 1
\right\}\right)>0,&&
\end{eqnarray}respectively.
If $p_{aa}=p_{bb}$, then (\ref{A3equal}) implies that both (\ref{eng}) and
(\ref{ing}) are  satisfied, and hence both (\ref{ab}) and
(\ref{ba}) hold. If $p_{aa}>p_{bb}$, then
(\ref{eng}), and subsequently \eqref{ab}, follow from (\ref{A3unequal-a}).
If $p_{aa}<p_{bb}$, then
(\ref{ing}), and subsequently \eqref{ba}, follow from (\ref{A3unequal-b}).
Hence,  at least one of the assumptions (\ref{ab}), (\ref{ba}) is always guaranteed to hold.

\underline{In  case 2}, (\ref{ab}) and (\ref{ba}) are equivalent to
\begin{eqnarray}\label{ahv}
P_a\left(\left\{x\in \mathcal{X}:f_a(x)p_{ba}>
f_b(x)p_{ab}\right)\right\}=P_a\left(\left\{x\in \mathcal{X}:{f_b(x)p_{ab}\over f_a(x)p_{ba}}< 1 \right\}\right)>0&&\\
\label{paavian}P_b\left(\left\{x\in \mathcal{X}:f_b(x)p_{ab}>
f_a(x)p_{ba}\right)\right\}=P_b\left(\left\{x\in \mathcal{X}:{f_a(x)p_{ba}\over f_b(x)p_{ab}}< 1 \right\}\right)>0,&&
\end{eqnarray} respectively.
Again, if $p_{aa}=p_{bb}$, then (\ref{ahv}) and (\ref{paavian})
both hold without further assumptions. If $p_{aa}>p_{bb}$, then
(\ref{ahv}) is automatically  satisfied. Likewise, (\ref{paavian}) holds if $p_{aa}<p_{bb}$.
Hence, one of the assumptions \eqref{ab}, \eqref{ba} is always guaranteed to hold.

\underline{In  case 3}, \eqref{ab} and \eqref{ba} write
\begin{align}\label{abmixture} & P_a\left(\left\{x\in \mathcal{X}:f_a(x)\pi_a>
f_b(x)\pi_b\right\}\right)>0,\\
\label{bamixture} & P_b\left(\left\{x \in \mathcal{X}: f_b(x)\pi_b>f_a(x)\pi_a\right\}\right)>0.
\end{align}
Assume $\pi_a\ge \pi_b$. Then, \eqref{A3equal} implies 
$\lambda\left(\left \{x\in \mathcal{X}: \pi_af_a(x)>\pi_bf_b(x)\right\}\right)>0$, which in 
turn implies \eqref{abmixture}.  
\end{proof}

Finally, we state and prove the main results for each of the three cases.
\\
\begin{figure}[htbp]
 \begin{center}
  \setlength{\unitlength}{0.00055in}
\begingroup\makeatletter\ifx\SetFigFont\undefined%
\gdef\SetFigFont#1#2#3#4#5{%
  \reset@font\fontsize{#1}{#2pt}%
  \fontfamily{#3}\fontseries{#4}\fontshape{#5}%
  \selectfont}%
\fi\endgroup%
{\renewcommand{\dashlinestretch}{30}
\begin{picture}(7808,6660)(0,-10)
\put(450,1088){\ellipse{150}{150}}
\put(4050,1088){\ellipse{150}{150}}
\put(5850,1088){\ellipse{150}{150}}
\put(7650,1088){\ellipse{150}{150}}
\put(1350,188){\ellipse{150}{150}}
\put(2250,1088){\ellipse{150}{150}}
\put(3150,1088){\ellipse{150}{150}}
\put(4950,188){\ellipse{150}{150}}
\put(6750,188){\ellipse{150}{150}}
\put(450,3788){\ellipse{150}{150}}
\put(450,2888){\ellipse{150}{150}}
\put(1350,3788){\ellipse{150}{150}}
\put(4050,2888){\ellipse{150}{150}}
\put(4950,3788){\ellipse{150}{150}}
\put(5850,2888){\ellipse{150}{150}}
\put(5850,3788){\ellipse{150}{150}}
\put(7650,2888){\ellipse{150}{150}}
\put(7650,3788){\ellipse{150}{150}}
\put(1350,2888){\ellipse{150}{150}}
\put(2250,3788){\ellipse{150}{150}}
\put(3150,3788){\ellipse{150}{150}}
\put(3150,2888){\ellipse{150}{150}}
\put(6750,3788){\ellipse{150}{150}}
\put(6750,2888){\ellipse{150}{150}}
\put(450,6488){\ellipse{150}{150}}
\put(450,5588){\ellipse{150}{150}}
\put(1350,6488){\ellipse{150}{150}}
\put(2250,5588){\ellipse{150}{150}}
\put(4050,5588){\ellipse{150}{150}}
\put(4950,6488){\ellipse{150}{150}}
\put(5850,5588){\ellipse{150}{150}}
\put(7650,5588){\ellipse{150}{150}}
\put(7650,6488){\ellipse{150}{150}}
\put(2250,6488){\ellipse{150}{150}}
\put(3150,6488){\ellipse{150}{150}}
\put(3150,5588){\ellipse{150}{150}}
\put(4950,5588){\ellipse{150}{150}}
\put(6750,6488){\ellipse{150}{150}}
\put(6750,5588){\ellipse{150}{150}}
\put(1350,5588){\ellipse{300}{300}}
\put(4050,6488){\ellipse{300}{300}}
\put(5850,6488){\ellipse{300}{300}}
\put(2250,2888){\ellipse{300}{300}}
\put(4050,3788){\ellipse{300}{300}}
\put(4950,2888){\ellipse{300}{300}}
\put(1350,1088){\ellipse{300}{300}}
\put(450,188){\ellipse{300}{300}} \put(2250,188){\ellipse{300}{300}}
\put(3150,188){\ellipse{300}{300}}
\put(4950,1088){\ellipse{300}{300}}
\put(5850,188){\ellipse{300}{300}}
\put(6750,1088){\ellipse{300}{300}}
\put(7650,188){\ellipse{300}{300}}
\put(4050,188){\ellipse{300}{300}} \drawline(4950,3788)(4950,3788)
\drawline(4950,6488)(4950,6488) \path(525,5588)(1200,5588)
\path(1500,5663)(2175,6413) \path(2325,6488)(3075,6488)
\path(3225,6488)(3900,6488) \path(4200,6488)(4875,6488)
\path(5025,6488)(5700,6488) \path(6000,6488)(6675,6488)
\path(6825,6488)(7575,6488) \dashline{60.000}(525,6488)(1275,6488)
\dashline{60.000}(1500,5588)(2175,5588)
\dashline{60.000}(2325,5588)(3075,5588)
\dashline{60.000}(3225,5588)(3975,5588)
\dashline{60.000}(4200,6338)(4875,5663)
\dashline{60.000}(5025,5588)(5775,5588)
\dashline{60.000}(6825,5588)(7575,5588)
\dashline{60.000}(6000,6338)(6675,5663)
\dashline{60.000}(525,3713)(1275,2963)
\dashline{60.000}(1425,2963)(2175,3713)
\dashline{60.000}(2400,3038)(3075,3713)
\dashline{60.000}(3225,3713)(3975,2963)
\dashline{60.000}(4200,3788)(4875,3788)
\dashline{60.000}(5100,2888)(5775,2888)
\dashline{60.000}(5925,2963)(6675,3713)
\dashline{60.000}(6825,3713)(7575,2963) \path(525,2963)(1275,3713)
\path(1425,3713)(2100,3038) \path(2400,2888)(3075,2888)
\path(3225,2963)(3900,3638) \path(4125,3713)(4875,2963)
\path(5100,3038)(5775,3713) \path(5925,3713)(6675,2963)
\path(6825,2963)(7575,3713) \path(600,338)(1200,938)
\path(1500,938)(2100,338) \path(2400,188)(3000,188)
\path(3300,188)(3900,188) \path(5025,1013)(5775,263)
\path(6000,338)(6675,1013) \path(6825,1013)(7575,263)
\dashline{60.000}(600,188)(1275,188)
\dashline{60.000}(1500,1088)(2175,1088)
\dashline{60.000}(2400,338)(3075,1013)
\dashline{60.000}(3300,338)(3975,1013)
\drawline(3975,1013)(3975,1013)
\dashline{60.000}(4200,188)(4875,188)
\dashline{60.000}(5100,1088)(5775,1088)
\dashline{60.000}(6000,188)(6675,188) \drawline(6675,188)(6675,188)
\dashline{60.000}(6900,1088)(7575,1088) \path(4125,263)(4800,938)
\put(0,5513){\makebox(0,0)[lb]{\smash{{{\SetFigFont{12}{14.4}{\rmdefault}{\mddefault}{\updefault}b}}}}}
\put(0,3788){\makebox(0,0)[lb]{\smash{{{\SetFigFont{12}{14.4}{\rmdefault}{\mddefault}{\updefault}a}}}}}
\put(0,2813){\makebox(0,0)[lb]{\smash{{{\SetFigFont{12}{14.4}{\rmdefault}{\mddefault}{\updefault}b}}}}}
\put(0,6488){\makebox(0,0)[lb]{\smash{{{\SetFigFont{12}{14.4}{\rmdefault}{\mddefault}{\updefault}a}}}}}
\put(0,1088){\makebox(0,0)[lb]{\smash{{{\SetFigFont{12}{14.4}{\rmdefault}{\mddefault}{\updefault}a}}}}}
\put(0,113){\makebox(0,0)[lb]{\smash{{{\SetFigFont{12}{14.4}{\rmdefault}{\mddefault}{\updefault}b}}}}}
\end{picture}
}
\end{center}
\caption{Distinct patterns of the Viterbi alignment in the two-state HMM:
Top:  Case 1, state can possibly change only at nodes (larger circles).
Middle: Case 2, states always alternate,  except possibly at 
nodes. Bottom: Case 3, every observation is a node.}
  \label{fig:2state}
\end{figure}
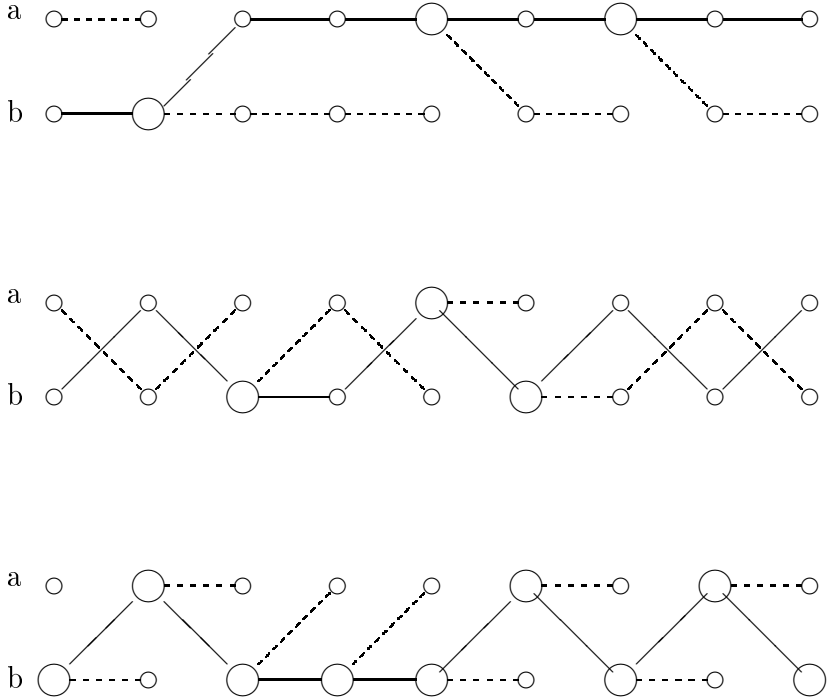
\subsection{Case 1}
First, note that condition (\ref{verystronglambda}) in this case is
equivalent to
\begin{equation}\label{scatman} \lambda\left(\left\{x\in \mathcal{X}:
p_{ba}f_a(x)p_{ab}>p_{bb}f_b(x)p_{bb}\right\}\right)>0,
\end{equation} As mentioned in \S\ref{sec:prelim}, condition \eqref{verystronglambda} need not hold in general.
Nonetheless, for the two-state HMM, we have the following Lemma.
\begin{lemma}\label{case1general}
In case 1, almost every realization of the two-state HMM has infinitely many strong
barriers.\end{lemma}
\begin{proof} Without loss of generality, assume $p_{aa}\geq p_{bb}$.
 Then (\ref{eng}) holds implying that
there exists $\epsilon>0$ such that
$$P_a({\cal X}_a)>0,
\quad \text{where}\quad  {\cal X}_a:=\left\{x\in {\cal X}: {f_b(x)p_{bb}\over
f_a(x)p_{aa}}< 1-\epsilon \right\}.$$ Let integer $k$ be sufficiently large for
$(1-\epsilon)^k<p_{ab}p_{ba}/(p_{aa}p_{bb})$ to hold.
Then every sequence $z_{1:k}\in {\cal X}_a^{k}$
satisfies
\begin{equation}\label{mambo}
\prod_{j=1}^k{f_b(z_j)p_{bb}\over f_a(z_j)p_{aa}}<
(1-\epsilon)^k<{p_{ab}p_{ba}\over p_{aa}p_{bb}}.
\end{equation}
Let $u> k$ be arbitrary and let $z_{0:k}\in {\cal X}_a^{k+1}$ be the last $k+1$ observations in a generic
sequence $x_{1:u}\in \mathcal{X}^{u-k-1}\times  {\cal X}_a^{k+1}$. To shorten the notation,
we write  $d_j(z_i)$ for  $\delta_{u-k+i}(j)$ for every
$i=0,1,\ldots,k$, $j=a,b$. Next, we show that $x_{u-k:u}$ contains at least one strong node, and
consequently, $z_{0:k}$ is a strong barrier.
Indeed, if none of the observations $x_{u-k:u}$ were a
strong $a$-node then we would have
$$d_b(z_k)=d_{b}(z_0)\prod_{j=1}^kf_b(z_j)p_{bb}.$$
Similarly, if none among the observations $x_{u-k+1:u}$ were a strong $b$-node, we
would have
$$\delta_u(a)\geq \delta_{u-k}(b)p_{ba}(\prod_{j=1}^kf_a(z_j))p_{aa}^{k-1}.$$
Hence,
$${\delta_u(b)\over \delta_u(a)}\leq
{\delta_{u-k}(b)p_{bb}(\prod_{j=1}^kf_b(z_j))p_{bb}^{k-1}\over
\delta_{u-k}(b)p_{ba}(\prod_{j=1}^k
f_a(z_j))p_{aa}^{k-1}}={\prod_{j=1}^k(f_b(z_j)p_{bb})\over
\prod_{j=1}^k(f_a(z_j)p_{aa})}{p_{aa}\over p_{ba}}$$ and by
(\ref{mambo})
$${\delta_u(b)\over \delta_u(a)}<{p_{ab}\over p_{bb}}$$
that contradicts  (\ref{yks}). Thus, at least one of $x_{u-k:u}$
must be a strong node. Since $P_{a}({\cal X}_a)>0$,
by ergodicity of HMM, almost every realization has infinitely many
barriers $z_{0:k}\in {\cal X}_a^{k+1}$, implying also that
every realization has infinitely many strong nodes.\end{proof}

The next Theorem refines the previous result.
\begin{theorem}\label{case1}Suppose the (transition matrix of the) two-state HMM meets
the condition of case 1.
If $p_{aa}\geq p_{bb}$, then almost every realization
has infinitely many strong $a$-barriers. (If $p_{aa}\leq p_{bb}$,
then almost every realization has infinitely many strong
$b$-barriers.)
\end{theorem}
\begin{proof} Let $p_{aa}\geq p_{bb}$ and use the notation of the proof of Lemma \ref{case1general}. 
First, we show that none of the observations $x_{k-u+1:u}$ is
a $b$-node. Indeed, since
$$d_b(z_1)=\max\{d_a(z_0)p_{ab},d_b(z_0)p_{bb}\}f_b(z_1),$$
 at least one of the following two
inequalities must hold:
\begin{equation}\label{eqn:x}
p_{ab}f_b(z_1)p_{ba}\geq p_{aa}f_a(z_1)p_{aa},\quad
p_{bb}f_b(z_1)p_{ba}\geq p_{ba}f_a(z_1)p_{aa}\end{equation} in order for $x_{u-k+1}$ to be a $b$-node.
However,  \eqref{eng} implies that $p_{ba}f_a(z_1)p_{aa}>p_{bb}f_b(z_1)p_{ba}$ and,
since $p_{bb}> p_{ab}$, we have
$p_{bb}f_b(z_1)p_{ba}>p_{ab}f_b(z_1)p_{ba}$. Hence, neither of the two inequalities \eqref{eqn:x}
holds.  Thus, $x_{u-k+1}$ cannot be a $b$-node,
and the same argument shows that none of the subsequent observations $x_{u-k+2},\ldots,x_u$
can be a $b$-node either.

The argument of the proof of Lemma \ref{case1general} then shows that
one of the observations in $x_{u-k:u}$ is a strong $a$-node and therefore 
$z_{0:k}$ is a strong $a$-barrier. The ergodic argument finishes the proof. (The same
argument with $a$ and $b$ swapped establishes the second part of the
Theorem.)\end{proof}

Note that the condition $p_{bb}\geq p_{aa}$ is sufficient but not
necessary for (\ref{ing}) to hold. In fact, for many 2-state HMMs, such as the one with additive
white Gaussian noise, both (\ref{eng}) and (\ref{ing}) hold for any
(positive) values of $p_{aa}$ and $p_{bb}$. On the other hand, it
might happen that one of the conditions \eqref{eng} and \eqref{ing}, say \eqref{ing}, fails.
This would mean $P_b\left(\left\{x \in \mathcal{X}: p_{bb}f_b(x)>p_{aa}f_a(x)\right\}\right)=0$ or, equivalently,
\begin{equation}\label{8}
\lambda \left(\left\{x \in \mathcal{X}:  p_{bb}f_b(x)>p_{aa}f_a(x)\right\}\right)=0.
\end{equation}
\begin{corollary}\label{cor} In case 1, equation \eqref{8} implies
that {\em almost every} sequence of observations has infinitely many
strong $a$-barriers and no strong $b$-nodes. Furthermore, equation \eqref{8} in case 1
implies that for {\em almost every} realization, if a $b$-node does occur, it occurs before the first $a$-node.
\end{corollary}
\begin{proof}
From the proof of Theorem \ref{case1}, it follows that no
observation $x\in \X$ such that $p_{bb}f_b(x)\leq p_{aa}f_a(x)$
(i.e. from the complement of the set in \eqref{8}) can be a strong $b$-node; a
closer inspection of the proof actually shows that even a weak (i.e. not strong) $b$-node cannot
occur after an $a$-node (since in case 1 $p_{bb}>p_{ba}$).  Theorem
\ref{case1} then implies that {\em almost every} sequence of observations has
infinitely many strong $a$-barriers.\end{proof}

Corollary \ref{cor} in its turn implies that starting with the first strong  $a$-node
onward, the Viterbi alignment $v(x_{1:n})$ stays in state $a$. As we have already mentioned,
Viterbi alignments need not be unique (see \citep{AVT4}), i.e. ties are possible in general,
and in this case, in particular, they are possible up until the first strong  $a$-node.
However, the impossibility of strong $b$-nodes in this case implies that the ties
can be broken in favor of $a$, resulting in the constant all $a$ alignment.

Theorem \ref{case1} is a generalization of Theorem 7 in \citep{caliebe2},
which basically states that in case 1, if (\ref{eng}) and (\ref{ing})
hold then under some additional assumptions (equal supports of $P_a$ and $P_b$ and further
conditions A2), almost every realization has infinitely many nodes. Thus,
\citep{caliebe2} stops short of realizing that in case 1
conditions (\ref{eng}) and (\ref{ing}) alone ensure the existence of
$a-$ and $b$-nodes. This results in \citep{caliebe2} invoking Theorem 2
of \citep{caliebe1} to prove the existence of nodes, hence
superfluous assumptions A1, A2. Also the proof
of Theorem 7 in \citep{caliebe1} could be simplified and shortened with
the help of the notions of nodes and barriers.
Finally, Corollary \ref{cor} generalizes Theorems 8 and 9 of \citep{caliebe2}.
\subsection{Case 2}Recall that we have been proving the existence of barriers without condition \eqref{verystronglambda}. Note that in case 2, condition \eqref{verystronglambda} becomes
\begin{equation*}
\lambda\left(\left\{x\in \mathcal{X}: p_{aa}f_a(x)p_{aa}>p_{ab}f_b(x)p_{ba}\right\}\right)>0.
\end{equation*}
Recall  (\S\ref{sec:prelim}) also that interchanging $a$ with $b$ gives a similar condition for strong
$b$-nodes to occur infinitely often in {\em almost every} realization.

It follows from \eqref{A3equal} that for some $\epsilon>0$, the sets
$${\cal X}_a:=\{x\in \mathcal{X}: f_a(x)(1-\epsilon)>f_b(x)\},\quad {\cal X}_b:=\{x\in \mathcal{X}:
f_a(x)<f_b(x)(1-\epsilon)\}$$ both have positive $\lambda$-measure.
Hence $P_a(\X_a)>0$ and $P_b(\X_b)>0$. Then, for $x_{1:2}\in {\cal X}_a\times {\cal X}_b$, the
following holds:
\begin{equation}\label{korrutis}
{f_b(x_1)f_a(x_2)\over f_a(x_1)f_b(x_2)}<(1-\epsilon)^2.
\end{equation}
\begin{lemma}\label{case2general}
In case 2, almost every realization has infinitely many strong
barriers.\end{lemma}
\begin{proof} Let ${\cal X}_a$ and ${\cal X}_b$ be as above. Choose $k$
sufficiently large for
$$(1-\epsilon)^{2k}<{p_{aa}p_{bb}\over p_{ba}p_{ab}}$$
to hold. Next, consider a sequence $z_{0:2k}\in \X^{2k+1}$, where $z_0,z_{2i}\in {\cal
X}_a$, $z_{2i-1}\in {\cal X}_b$, for every $i=1,\ldots, k$. We show that for every $u>2k$, 
every sequence of observations $x_{1:u}\in \mathcal{X}^u$  such
that $x_{u-2k:u}=z_{0:2k}$, contains a strong node, making  $z_{0:2k}$  a strong barrier.

The choice of $k$ and $z_{0:2k}$ implies
\begin{equation}\label{uhhuu}
{\prod_{i=1}^kp_{ba}f_a(z_{2i-1})p_{ab}f_b(z_{2i})\over
\prod_{i=1}^kp_{ab}f_b(z_{2i-1})p_{ba}f_a(z_{2i})}<(1-\epsilon)^{2k}<{p_{bb}p_{aa}\over
p_{ba}p_{ab}}.
\end{equation}
If there is no strong node among $x_{u-2k:u}$, then
$$d_b(z_{2k})=d_b(z_0)\prod_{i=1}^kp_{ba}f_a(z_{2i-1})p_{ab}f_b(z_{2i})$$
and
$$d_a(z_{2k})\geq d_b(z_0){p_{bb}\over p_{ab}}\prod_{i=1}^kp_{ab}f_b(z_{2i-1}) p_{ba}f_a(z_{2i}).$$
Hence, by (\ref{uhhuu})
$${d_b(z_{2k})\over d_a(z_{2k})}\leq {\prod_{i=1}^kp_{ba}f_a(z_{2i-1})p_{ab}f_b(z_{2i})\over
{p_{bb}\over p_{ab}}\prod_{i=1}^kp_{ab}f_b(z_{2i-1})
p_{ba}f_a(z_{2i})}<{p_{aa}\over p_{ba}}$$ which contradicts
(\ref{kaks}). \end{proof}

Next, we refine this result. Without loss of generality assume
$p_{ba}\geq p_{ab}$.  Therefore
\begin{equation}\label{vorratused}
p_{ab}p_{aa}\geq p_{ba}p_{bb},
\end{equation}
and also, for every $x\in \X_a$,
\begin{equation}\label{vorratus-x}
p_{ba}f_a(x)> p_{ab}f_b(x).
\end{equation}
Hence, \eqref{ahv} holds. We multiply the right side of
\eqref{vorratus-x} by $p_{ba}p_{bb}$ and the left side by
$p_{ab}p_{aa}$,  and use \eqref{vorratused} to obtain 
\begin{equation}\label{essa}
f_a(x)p_{aa}>f_b(x)p_{bb}.
\end{equation}
Finally, for $x\in {\cal X}_b$, we have
\begin{equation}\label{vorratus-y}
f_a(x)< f_b(x).
\end{equation}
We 
will need the following Lemma.
\begin{lemma}\label{lemma:case2:aux}
Assume (in addition to being in case 2) that $p_{ab}\le p_{ba}$.
\begin{description}
  \item[a)] In any pair of observations $z_{1:2}\in
    \mathcal{X}_a\times \mathcal{X}_b$, $z_1$ is not a $b$-node.
  \item[b)] In any pair of observations $z_{2:3}\in
    \mathcal{X}_b\times \mathcal{X}_a$, if $z_2$ is a $b$-node,
then $z_3$ is a strong $a$-node.
\end{description}
\end{lemma}
\begin{proof}
Assume that $p_{ab}\le p_{ba}$, and consider ${\bf a)}$. 
First note that since we are in case 2,
$z_1$ is  a $b$-node if and only if
\begin{equation}\label{b-node}
d_b(z_1)p_{bb}\geq d_a(z_1)p_{ab}.
\end{equation}
Suppose first that $z_0$ is not a
node, in which case $d_b(z_1)=d_a(z_0)p_{ab}f_{b}(z_{1})$ and
$d_a(z_1)=d_b(z_0)p_{ba}f_{a}(z_{1})$. Then
\begin{align*}
d_a(z_1)p_{ab}&=d_b(z_0)p_{ba}f_{a}(z_{1})p_{ab}\geq
d_a(z_0)p_{aa}f_{a}(z_{1})p_{ab}\\&>d_a(z_0)p_{bb}f_{b}(z_{1})p_{ab}=
d_a(z_0)p_{ab}f_{b}(z_{1})p_{bb}=d_b(z_1)p_{bb}.
\end{align*}
The first inequality above follows from the recursion property
\eqref{eqn:recurse} 
of scores $\delta$, whereas
the second one follows from \eqref{essa}. Thus, when $z_0$ is
not a node, $z_1$ cannot be a $b$-node. Similarly, supposing
that $z_0$ is an $a$-node, we obtain that $z_1$ is not  a $b$-node.  
Suppose finally that $z_0$ is a $b$-node. Then
$d_b(z_1)=d_b(z_0)p_{bb}f_{b}(z_{1})$ and
$d_a(z_1)=d_b(z_0)p_{ba}f_{a}(z_{1})$. Applying
consecutively $p_{bb}<p_{ab}$, (\ref{vorratus-x}) and
$p_{bb}<p_{ab}$ again, we obtain:
$p_{bb}f_{b}(z_{1})p_{bb}<p_{ab}f_{b}(z_{1})p_{bb}\leq
p_{ba}f_{a}(z_{1})p_{bb}<p_{ba}f_{a}(z_{1})p_{ab}.$ Thus, contrary to \eqref{b-node}
$$d_b(z_1)p_{bb}=d_b(z_0)p_{bb}f_{b}(z_{1})p_{bb}<d_b(z_0)p_{ba}f_{a}(z_{1})p_{ab}
=d_a(z_1)p_{ab},$$
that is, $z_1$ is not a $b$-node in this case either.
Let us now prove {\bf b)}. If $z_{2}$ is a $b$-node, then
$d_a(z_3)=d_b(z_2)p_{ba}f_a(z_{3})$ and
$d_b(z_3)=d_b(z_2)p_{bb}f_b(z_{3})$. By (\ref{essa}), we
now have $d_a(z_3)p_{aa}=
d_b(z_2)p_{ba}f_a(z_3)p_{aa}>d_b(z_2)p_{bb}f_b(z_3)p_{ba}=d_b(z_3)p_{ba}.$
Similarly to the argument regarding $b$-nodes guaranteed by 
\eqref{b-node} above,  we now have $d_a(z_3)>d_b(z_3)$, implying
$d_a(z_3)p_{ab}>d_b(z_3)p_{bb}$. Thus $z_3$ is a strong $a$-node.
\end{proof}
\begin{theorem}\label{case2} If $p_{ba}\geq p_{ab}$, then almost every
realization has infinitely many strong $a$-nodes. If $p_{ba}\leq
p_{ab}$, then almost every realization has infinitely many strong
$b$-nodes.\end{theorem}
\begin{proof} Assume again that $p_{ba}\geq p_{ab}$. 
Let $z_{0:2k}$ be as in the proof of  Lemma \ref{case2general}
and attach one more element $z_{2k+1}\in \X_b$ to the end. 
Thus, $z_{2i}\in \X_a$ and $z_{2i+1}\in \X_b$, $i=0,1,\ldots,k$. 

From (the proof of) Lemma \ref{case2general} we know that $z_{0:2k}$
contains at least one strong node. If this is an $a$-node,
then the theorem is proven. Otherwise this is a $b$-node,
which, according to part {\bf a)} of
Lemma~\eqref{lemma:case2:aux}, can only be among $z_1$, $z_3$, \ldots, 
$z_{2k-1}$.  Applying  part {\bf b)} of 
Lemma~\eqref{lemma:case2:aux} shows that  there must also be a strong
$a$-node $z_2$, $z_4$, \ldots, $z_{2k}$.  
Invoking ergodicity again finishes the proof. 

Clearly, swapping $a$ and $b$ in the above discussion following
the proof of Lemma~\ref{case2general}, establishes the other part of
the theorem.\end{proof}

Inequality (\ref{vorratused}) guarantees \eqref{ahv}.
Often, the model is such that in addition to \eqref{ahv},
\eqref{paavian} also holds. However, to apply the previous proof
(i.e. of Theorem~\ref{case2}) to
guarantee  the simultaneous existence of infinitely many strong $a$
and $b$-nodes, we would need the following counterpart of \eqref{essa}:
$P_b (\{x\in \X: f_b(x)p_{ab}>f_a(x)p_{ba},
f_b(x)p_{bb}>f_a(x)p_{aa}\})>0$, 
which is stronger than \eqref{paavian}.
However, this previous condition is indeed often met, 
resulting in infinitely many strong $a$- and $b$-nodes
(in {\em almost every} realization $x_{1:\infty}$).

Lemma \ref{case2general} appears without proof as Theorem 10 in
\citep{caliebe2}. 
The author of \citep{caliebe2} actually suggests that Theorem 10 and
 other results for case 2 are analogous to the corresponding results
 for case 1, mainly Theorem 7 (of the same work). It is further 
stated in \citep{caliebe2} that the proofs of those results are not
given as they ``are very similar'' to the corresponding proofs in case 1.
Our present workings actually show that case 2
is quite  dissimilar to case 1 (due to the fluctuating  nature of 
the typical Viterbi alignment) and in particular requires a more 
careful treatment. Note that, even if Theorem 10 in \citep{caliebe2} 
assumed \eqref{ab} and \eqref{ba} (as  Theorem 7 in \citep{caliebe2}
does) to help one to prove this Theorem by analogy to Theorem 7, 
it is still not clear how the two proofs could be very similar. 
\subsubsection{Case 3 (the mixture model) }
Recall that every observation in this case is a (not necessarily strong) node.
Furthermore, every observation from $\{x\in \mathcal{X}: \pi_af_a(x)>\pi_bf_b(x)\}$ is a strong
$a$ node. Thus,  we have the following counterpart of Theorems \ref{case1} and
\ref{case2}.
\begin{theorem}\label{case3} If $\pi_{a}\geq \pi_{b}$, then almost every
realization has infinitely many strong $a$-nodes. If $\pi_{a}\leq
\pi_{b}$, then almost every realization has infinitely many strong
$b$-nodes.\end{theorem}
\section{Conclusion}
In summary, we have proved Theorem~\ref{maintheorem} stated below and
providing a basis for the piecewise construction and asymptotic analysis
of the Viterbi alignments of two-state HMMs.
\begin{theorem}\label{maintheorem}
Almost every realization of the two-state HMM has infinitely many strong
barriers.   Furthermore
\begin{enumerate}[a)]
\item if the transition probabilities satisfy $p_{aa}\ge p_{ba}$  then
(almost every realization of) the chain has infinitely many strong $s$-barriers
where $s$ is such that $p_{ss}=\max\{p_{aa},p_{bb}\}$,
\item otherwise (i.e. if $p_{aa}<p_{ba}$)
(almost every realization of) the chain has infinitely many strong $s$-barriers
where $s$ is such that $p_{ts}=\max\{p_{ab},p_{ba}\}$ (for some $t\in S$).
\end{enumerate}
\end{theorem}


\begin{thebibliography}{29}
\expandafter\ifx\csname natexlab\endcsname\relax\def\natexlab#1{#1}\fi
\expandafter\ifx\csname url\endcsname\relax
  \def\url#1{\texttt{#1}}\fi
\expandafter\ifx\csname urlprefix\endcsname\relax\def\urlprefix{URL }\fi

\bibitem[{Baum and Petrie(1966)}]{baumhmm}
Baum, L.~E., Petrie, T., 1966. Statistical inference for probabilistic
  functions of finite state {M}arkov chains. Ann. Math. Statist. 37,
  1554--1563.

\bibitem[{Bilmes(1998)}]{emtutorial}
Bilmes, J., 1998. A gentle tutorial of the {EM} algorithm and its application
  to parameter estimation for {G}aussian mixture and hidden {M}arkov models.
  Tech. Rep. 97--021, International Computer Science Institute, Berkeley, CA,
  USA.

\bibitem[{Caliebe(2006)}]{caliebe2}
Caliebe, A., 2006. Properties of the maximum a posteriori path estimator in
  hidden {M}arkov models. IEEE Trans. Inform. Theory 52~(1), 41--51.

\bibitem[{Caliebe and R{\"o}sler(2002)}]{caliebe1}
Caliebe, A., R{\"o}sler, U., 2002. Convergence of the maximum a posteriori path
  estimator in hidden {M}arkov models. IEEE Trans. Inform. Theory 48~(7),
  1750--1758.

\bibitem[{Capp{\'e} et~al.(2005)Capp{\'e}, Moulines, and Ryd{\'e}n}]{HMM}
Capp{\'e}, O., Moulines, E., Ryd{\'e}n, T., 2005. Inference in hidden {M}arkov
  models. Springer Series in Statistics. Springer, New York, with Randal Douc's
  contributions to Chapter 9 and Christian P. Robert's to Chapters 6, 7 and 13,
  With Chapter 14 by Gersende Fort, Philippe Soulier and Moulines, and Chapter
  15 by St\'ephane Boucheron and Elisabeth Gassiat.

\bibitem[{Durbin et~al.(1998)Durbin, Eddy, A., and Mitchison}]{BioHMM2}
Durbin, R., Eddy, S., A., K., Mitchison, G., 1998. Biological {S}equence
  {A}nalysis: {P}robabilistic {M}odels of {P}roteins and {N}ucleic {A}cids.
  Cambridge University Press.

\bibitem[{Eddy(2004)}]{natHMM}
Eddy, S., 2004. What is a hidden {M}arkov model? Nature Biotechnology 22~(10),
  1315 -- 1316.

\bibitem[{Ephraim and Merhav(2002)}]{HMP}
Ephraim, Y., Merhav, N., 2002. Hidden {M}arkov processes. IEEE Trans. Inform.
  Theory 48~(6), 1518--1569, special issue on Shannon theory: perspective,
  trends, and applications.

\bibitem[{Genon-Catalot et~al.(2000)Genon-Catalot, Jeantheau, and
  Lar{\'e}do}]{volatility}
Genon-Catalot, V., Jeantheau, T., Lar{\'e}do, C., 2000. Stochastic volatility
  models as hidden {M}arkov models and statistical applications. Bernoulli
  6~(6), 1051--1079.

\bibitem[{Huang et~al.(1990)Huang, Ariki, and Jack}]{vanaraamat}
Huang, X., Ariki, Y., Jack, M., 1990. Hidden {M}arkov models for speech
  recognition. Edinburgh University Press, Edinburgh, UK.

\bibitem[{Jelinek(1976)}]{jelinek0}
Jelinek, F., 1976. Continuous speech recognition by statistical methods. Proc.
  IEEE 64, 532--556.

\bibitem[{Jelinek(2001)}]{jelinek}
Jelinek, F., 2001. Statistical methods for speech recognition. The MIT Press,
  Cambridge, MA, USA.

\bibitem[{Ji and Bilmes(2006)}]{bilmes2}
Ji, G., Bilmes, J., 2006. Backoff model training using partially observed data:
  {A}pplication to dialog act tagging. In: Proc. Human Language Techn. Conf.
  NAACL, Main Conference. Association for Computational Linguistics, New York
  City, USA, pp. 280--287.
\newline\urlprefix\url{http://www.aclweb.org/anthology/N/N06/N06-1036}

\bibitem[{Juang and Rabiner(1990)}]{Rabiner90}
Juang, B.-H., Rabiner, L., 1990. The segmental {K}-means algorithm for
  estimating parameters of hidden {M}arkov models. IEEE Trans. Acoust. Speech
  Signal Proc. 38~(9), 1639--1641.

\bibitem[{Koloydenko et~al.(2007)Koloydenko, K{\"a}{\"a}rik, and Lember}]{AVT3}
Koloydenko, A., K{\"a}{\"a}rik, M., Lember, J., 2007.
  \href{http://www.springerlink.com/content/t5882055g807t574/}{On adjusted
  {V}iterbi training}. Acta Appl. Math. 96~(1-3), 309--326.

\bibitem[{Krogh(1998)}]{BioHMM1}
Krogh, A., 1998. Computational {M}ethods in {M}olecular {B}iology. Elsevier
  Science, Ch. An {I}ntroduction to {H}idden {M}arkov {M}odels for {B}iological
  {S}equences.

\bibitem[{Lember and Koloydenko(2007)}]{AVT1}
Lember, J., Koloydenko, A., 2007.
  \href{http://journals.cambridge.org/action/displayAbstract?fromPage=online&a%
id=1293756&fulltextType=RA&fileId=S0269964807000083} {Adjusted {V}iterbi
  training: {A} proof of concept}. Probab. Eng. Inf. Sci. 21~(3), 451--475.

\bibitem[{Lember and Koloydenko(2008)}]{AVT4}
Lember, J., Koloydenko, A., 2008.
  \href{http://projecteuclid.org/handle/euclid.bj/1202492790/} {The {A}djusted
  {V}iterbi training for hidden {M}arkov models}. Bernoulli 14~(1), 180--206.

\bibitem[{Leroux(1992)}]{leroux}
Leroux, B.~G., 1992. Maximum-likelihood estimation for hidden {M}arkov models.
  Stochastic Process. Appl. 40~(1), 127--143.

\bibitem[{Lomsadze et~al.(2005)Lomsadze, Ter-Hovhannisyan, Chernoff, and
  Borodovsky}]{findgene}
Lomsadze, A., Ter-Hovhannisyan, V., Chernoff, V., Borodovsky, M., 2005. Gene
  identification in novel eukaryotic genomes by self-training algorithm.
  Nucleic Acids Res. 33~(20), 6494--6506.

\bibitem[{McDermott and Hazen(2004)}]{mitEM2}
McDermott, E., Hazen, T., 2004. Minimum classification error training of
  landmark models for real-time continuous speech recognition. In: Proc.
  ICASSP.

\bibitem[{Ney et~al.(1994)Ney, Steinbiss, Haeb-Umbach, and {\em
  et.~al.}}]{philips}
Ney, H., Steinbiss, V., Haeb-Umbach, R., {\em et.~al.}, 1994. An overview of
  the {P}hilips research system for large vocabulary continuous speech
  recognition. Int. J. Pattern Recognit. Artif. Intell. 8~(1), 33--70.

\bibitem[{Och and Ney(2000)}]{ochney}
Och, F., Ney, H., 2000. Improved statistical alignment models. In: Proc. 38th
  Ann. Meet. Assoc. Comput. Linguist. Assoc. Comput. Linguist., pp. 440 -- 447.

\bibitem[{Padmanabhan and Picheny(2002)}]{IBM2002}
Padmanabhan, M., Picheny, M., 2002. Large-vocabulary speech recognition
  algorithms. Computer 35~(4), 42 -- 50.

\bibitem[{Rabiner and Juang(1993)}]{raamat}
Rabiner, L., Juang, B., 1993. Fundamentals of speech recognition.
  Prentice-Hall, Inc., Upper Saddle River, NJ, USA.

\bibitem[{Rabiner et~al.(1986)Rabiner, Wilpon, and Juang}]{Rabiner86}
Rabiner, L., Wilpon, J., Juang, B., 1986. A segmental {K}-means training
  procedure for connected word recognition. AT\&T {T}ech. J. 64~(3), 21--40.

\bibitem[{Shu et~al.(2003)Shu, Hetherington, and Glass}]{mitEM1}
Shu, I., Hetherington, L., Glass, J., 2003. Baum-{W}elch training for
  segment-based speech recognition. In: Proc. IEEE ASRU Workshop. St. Thomas,
  U. S. Virgin Islands,
  {http://groups.csail.mit.edu/sls/publications/2003/ASRU\_Shu.pdf}, pp.
  43--48.

\bibitem[{Steinbiss et~al.(1995)Steinbiss, Ney, Aubert, and {\em
  et.~al.}}]{philips2}
Steinbiss, V., Ney, H., Aubert, X., {\em et.~al.}, 1995. The {P}hilips research
  system for continuous-speech recognition. Philips J. Res. 49, 317--352.

\bibitem[{Str\"{o}m et~al.(1999)Str\"{o}m, Hetherington, Hazen, Sandness, and
  Glass}]{strom}
Str\"{o}m, N., Hetherington, L., Hazen, T., Sandness, E., Glass, J., 1999.
  Acoustic modeling improvements in a segment-based speech recognizer. In:
  Proc. IEEE ASRU Workshop. pp. 139--142.

\end{thebibliography}
\end{document}